\documentclass[reqno]{amsart}

\usepackage{yhmath}
\usepackage{amsmath}
\usepackage{amsfonts}
\usepackage{amsthm}
\usepackage{amssymb}
\usepackage{eucal}

\newtheorem{Theorem}{Theorem}[section]

\theoremstyle{definition}
\newtheorem{Remark}[Theorem]{Remark}

\numberwithin{equation}{section}

\newcommand{\R}{{\mathbb{R}}}
\newcommand{\N}{\mathbb{N}}

\newcommand{\fP}{{\bf P}}

\newcommand{\cP}{{\mathcal P}}
\newcommand{\cQ}{{\mathcal Q}}
\newcommand{\cL}{{\mathcal L}}
\newcommand{\cF}{{\mathcal F}}
\newcommand{\cK}{{\mathcal K}}

\newcommand{\pt}{\tilde{p}}

\newcommand{\PP}{\mathbb{P}}       
\newcommand{\PPt}{\tilde{\PP}}
\newcommand{\G}{{\Gamma}}

\newcommand{\At}{{\tilde{A}}}                                    
\newcommand{\Bt}{{\tilde{B}}}                      

\newcommand{\rank}{\mathrm{rank}}
\newcommand{\Id}{\mathrm{Id}}

\begin{document}

\title{Two variable deformations of the Chebyshev measure}

\author[J.~Geronimo]{Jeffrey~S.~Geronimo}\thanks{JSG was partially
supported by an NSF grant}
\address{School of Mathematics, Georgia Institute of Technology,
Atlanta, GA 30332--0160, USA}
\email{geronimo@math.gatech.edu}

\author[P.~Iliev]{Plamen~Iliev}
\address{School of Mathematics, Georgia Institute of Technology,
Atlanta, GA 30332--0160, USA}
\email{iliev@math.gatech.edu}

\date{December 16, 2006}

\begin{abstract} We construct one and two parameter deformations of the two 
dimensional Chebyshev polynomials with simple recurrence coefficients, 
following the algorithm in \cite{DGIM}. Using inverse scattering techniques, 
we compute the corresponding orthogonality measures.

\end{abstract}

\maketitle

\section{Introduction}
The general theory of bivariate orthogonal polynomials goes back 
to the work of Jackson \cite{J}. Special examples have
arisen in studies related to symmetric groups \cite{D, Ko3,Mac}), as
extensions of one variable polynomials \cite{FPP,Ko2} and as eigenfunctions of
partial differential equations \cite{KS,Ko1,KKL, KLL}. An updated account of 
the theory can be found in the books \cite{DX,S}.

A major difficulty encountered in the theory of orthogonal polynomials of more 
than one variable is which monomial ordering to use. Starting with \cite{J}, 
the preferred  ordering is the total degree ordering and for polynomials with 
the same total degree the ordering is lexicographical that is 
$$(k,l)<_{\text{totaldeg}}(k_1,l_1)$$ 
if 
$$k+l<k_1+l_1\ {\rm or}\ (k+l=k_1+l_1\ (k,l)<_{\text{lex}} (k_1,l_1)),$$
where lex means the lexicographical ordering (see below). In other words, we 
apply the Gram-Schmidt process to the polynomials ordered as follows
$\{1,y,x,y^2,xy,x^2,\dots\}$. If we denote by $\fP_n(x,y)$ the $(n+1)$ 
dimensional vector with components the orthonormal polynomials of total degree 
$n$, then the multiplications by $x$ and $y$ are given by three term 
recurrence relations
\begin{align}
&x\fP_n=A_{x,n}\fP_{n+1}+B_{x,n}\fP_n+A^t_{x,n-1}\fP_{n-1}\label{1.1}\\
&y\fP_n=A_{y,n}\fP_{n+1}+B_{y,n}\fP_n+A^t_{y,n-1}\fP_{n-1},\label{1.2}
\end{align}
where $A_{x,n}$, $A_{y,n}$ are $(n+1)\times (n+2)$ matrices such that 
$\rank(A_{x,n})=\rank(A_{y,n})=n+1$ and $B_{x,n}$, 
$B_{y,n}$ are symmetric $(n+1)\times (n+1)$ matrices. Notice that the Jacobi 
matrices corresponding to multiplications by $x$ and $y$ will commute which 
amounts to certain commutativity relations between the matrices defined above, 
see \cite{DX} for more details.

Inspired by the progress made in \cite{GW1} in connection with 
Fej\'er-Riesz factorizations and autoregressive filters in two variables, 
an alternative way to approach two dimensional orthogonal polynomials 
was proposed in \cite{DGIM} by relating them to the theory of matrix valued 
orthogonal polynomials. This can be accomplished by using the lexicographical 
ordering,
$$
(k,l)<_{\rm lex} (k_1,l_1)\Leftrightarrow k<k_1\mbox{ or }
(k=k_1\mbox{ and } l<l_1),
$$ or the reverse lexicographical ordering
$$
(k,l)<_{\rm revlex} (k_1,l_1)\Leftrightarrow
(l,k)<_{\rm lex} (l_1,k_1),
$$
to arrange 
the monomials. This naturally connects the theory of bivariate orthogonal 
polynomials to doubly Hankel matrices. In particular, this led to an 
alternate parametrization for positive doubly Hankel matrices and an 
algorithm for constructing their inverses. 

In the present paper we apply the algorithm from \cite{DGIM} to construct 
deformations of the two dimensional Chebyshev polynomials with relatively 
simple recurrence coefficients. These polynomials can be thought of 
as two dimensional analogs of the  Bernstein-Szeg\"o polynomials. 
Using the connection of the above theory with matrix valued orthogonal 
polynomials and the inverse scattering techniques developed in \cite{Ger} 
we show how one can obtain the orthogonality measure from the recurrence 
coefficients.

The paper is organized as follows. In the next section we recall the main 
ingredients needed for the construction. In particular, we review the 
parametrization of the doubly Hankel matrices in \cite{DGIM} and the 
algorithm which gives the recurrence coefficients. We also sketch 
the inverse scattering techniques in the matrix case  and 
the connection with the Darboux transformation. 
In Sections 3 and 4, respectively, we present one parameter and two parameter 
deformations of the recurrence coefficients for the two dimensional Chebyshev 
polynomials and we derive the orthogonality measure.

\section{Bivariate Orthogonal Polynomials}                    %
Let $\mu(x,y)$ be a Borel measure supported on $\R^2$, such that 
$$\int_{\R^2}f(x,y)d\mu(x,y)<\infty$$ for every polynomial $f(x,y)$.
 For every nonnegative integer $m$ we order the monomials $x^i y^j$, 
$0\le i\le n,\ 0\le j\le m$
lexicographically i.e.
$$\{1,y,\dots, y^{m},x,xy,\dots,xy^{m},x^2,\dots\}.$$ 
We construct the moment matrix
\begin{equation}\label{hankelone}
H_{n,m}=
\left[
\begin{matrix}
H_0 & H_1 & \adots & H_n
\\
H_1 & H_2 &  & H_{n+1}
\\
\adots &  & \adots & \adots \cr H_n & H_{n+1} & \adots & H_{2n}
\end{matrix}
\right],
\end{equation}
where each $H_i$ is an $(m+1)\times(m+1)$ matrix of the form
\begin{equation}\label{hankeltwo}
H_i=
\left[
\begin{matrix}
h_{i,0}& h_{i1} & \adots & h_{i,m}\\
h_{i1} & h_{i2}& \adots & \\ 
\adots & \adots & & \adots \\ 
h_{i,m} &  & \adots & h_{i,2m}
\end{matrix}
\right],\qquad i=0,\dots, 2n.
\end{equation}
Thus $H_{n,m}$ is a block Hankel matrix where each block is a
Hankel matrix so it has a doubly Hankel structure. If the reverse
lexicographical ordering is used in place of the lexicographical
ordering we obtain another moment matrix $\tilde H_{n,m}$ where
the roles of $n$ and $m$ are interchanged. We will assume that $H_{n,m}$ is 
positive definite for all $n$ and $m$. For every nonnegative integer $m$ we 
apply the Gram-Schmidt process to the basis of monomials ordered as above 
and define the orthonormal polynomials $p_{n,m}^l(x,y),\
0\le n,\, 0\le m , \, 0\le l\le m,$ by the equations,
\begin{equation}\label{sorthogonal}
\begin{split}
&\int_{\R^2} p_{n,m}^l x^i y^jd\mu(x,y)=0, \quad  0\le i<n\ {\rm and }\
0\le j\le m\ {\rm\ or}\ i=n\ {\rm and }\ 0\le j< l,
\\
&\int_{\R^2} p_{n,m}^l p_{n,m}^l d\mu(x,y)=1,
\end{split}
\end{equation}
and
\begin{equation}\label{sorthogonaldeg}
p_{n,m}^{l}(x,y) = k^{n,l}_{n,m,l} x^n y^l + \sum_{(i,j)<_{\rm
lex}(n,l)} k^{i,j}_{n,m,l}x^iy^j.
\end{equation}
With the convention $k^{n,l}_{n,m,l}>0$, the above equations uniquely
specify $p^l_{n,m}$. Polynomials orthonormal with respect to
$d\mu(x,y)$ but using the reverse lexicographical ordering will be
denoted by $\tilde p^l_{n,m}$. They are uniquely determined by the
above relations with the roles of $n$ and $m$ interchanged. Set,
\begin{equation}\label{vectpoly}
\PP_{n,m}=\left[\begin{matrix} p_{n,m}^{0}\\ p_{n,m}^{1}\\
\vdots\\ p_{n,m}^{m} \end{matrix}\right].
\end{equation}

The $\PP_{n,m}$ may be obtained in an alternate manner as follows. We 
associate an $(m+1)\times(m+1)$ matrix 
valued measure $dM^{m+1}(x)$ by taking 
\begin{equation}\label{matrixmeasure}
dM^{m+1}(x)=\int_{\R}[1,y,\dots,y^{m}]^td\mu_y(x,y)[1,y,\dots,y^{m}],
\end{equation}
where the above integral is with respect to $y$. Let us denote by 
$\{P_{n,m}\}_{n=0}^{\infty}$ the sequence of $(m+1)\times(m+1)$ matrix valued 
polynomials satisfying
\begin{align}\label{defp}
P_{n,m}(x)=K_{n,n}^m x^n + {\rm lower\ order\ terms},\\\
\int_{\R}P_{n,m}(x)dM^{m+1}(P_{k,m}(x))^t=\delta_{k,m}I_{m+1},
\end{align} 
with $K_{n,n}^m$ a lower triangular matrix with strictly positive  
diagonal entries. The above conditions uniquely specify these left matrix 
valued orthogonal polynomials and it follows that 
\begin{equation}\label{connectionform}
\PP_{n,m}(x,y)=P_{n,m}(x)[1,y,\dots,y^m]^t.
\end{equation} 
From equation~\eqref{matrixmeasure} we see that $M^{m+1}$ is a Hankel matrix 
and we have,
\begin{Theorem}\label{matrix}
Let $M^{m+1}(x)$ be a matrix measure supported on some interval 
$I\subset\R$ then $M^{m+1}(x)$ is a Hankel measure if and only if
there exists a bounded positive  linear functional 
$\cF_{2m}:C(I)\times\{1,\ldots,y^{2m}\}\to\R$ with 
$(\int_I x^k d M^{m+1})_{i,j}=\cF_{2m}(x^k y^{i+j})$.
\end{Theorem}

\begin{proof}
If we begin with $M^{m+1}(x)$ set $H_i=\int_I x^i dM^{m+1}$ and let 
$H_{n,m}$ be the $(n+1)(m+1)\times(n+1)(m+1)$ matrix given by 
equation~\eqref{hankelone}. Then $H_{n,m}$ is a positive definite doubly 
Hankel matrix which implies via the above formula connecting $M^{m+1}$ and 
$\cF_{2m}$ and the density of the polynomials in $C(I)$ that $\cF_{2m}$ is a 
bounded positive  linear functional which is uniquely defined. The converse 
follows from the relation between $\cF_{2m}$ and $M^{m+1}$ given above.
\end{proof}

The following Theorem was proved in \cite{DGIM}, 
\begin{Theorem}\label{recurrence} 
Given $\{\PP_{n,m}\}$ and 
$\{\tilde\PP_{n,m}\}$ the following recurrence formulas hold, 
\begin{align}
& x \PP_{n,m} = A_{n+1,m} \PP_{n+1,m} + B_{n,m} \PP_{n,m} +A_{n,m}^\top 
\PP_{n-1,m}, 
                     \label{2.10}\\
&\G_{n,m} \PP_{n,m} = \PP_{n,m-1} - \cK_{n,m}\PPt_{n-1,m},
\label{2.11}
\\
& J_{n,m}^1 \PP_{n,m} = y \PP_{n,m-1} + J_{n,m}^2
\PPt_{n-1,m} + J_{n,m}^3 \PPt_{n-1,m-1},
\label{2.12}
\\
&\PP_{n,m}=I_{n,m} \PPt_{n,m} + \G^{\top}_{n,m} \PP_{n,m-1}
\label{2.13}
\end{align}

with
\begin{align*}
A_{n,m} & = \langle x\PP_{n-1,m}, \PP_{n,m}\rangle, 
\\
B_{n,m} & = \langle x\PP_{n,m}, \PP_{n,m}\rangle 
\\
J^1_{n,m} & = \langle y \PP_{n,m-1}, \PP_{n,m} \rangle 
\\
J^2_{n,m} & = -\langle y\PP_{n,m-1}, \PPt_{n-1,m}\rangle 
\\
J^3_{n,m} & = -\langle y \PP_{n,m-1}, \PPt_{n-1,m-1}\rangle 
\\
\G_{n,m} & = \langle \PP_{n,m-1}, \PP_{n,m} \rangle 
\\
\cK_{n,m} & = \langle \PP_{n,m-1}, \PPt_{n-1,m}\rangle 
\\
I_{n,m} & = \langle \PP_{n,m}, \PPt_{n,m} \rangle 
\end{align*}
where for every two vector valued polynomials $f, g$ we set
\begin{equation}\label{inner}
\langle f,g\rangle =\int_{\R^2} f(x,y) g^t(x,y) d\mu(x,y).
\end{equation}
Similar formulas hold for $\PPt_{n,m}(x,y)$.

\end{Theorem}

Equation~\eqref{2.10} comes from the fact that the vector polynomials can be 
related to matrix orthogonal polynomials as discussed  above. Given the 
coefficients in this equation it is possible to compute the vector polynomials 
in a strip of size $m+1$ in the $n$ direction. Its tilde counterpart obtained 
from the reverse lexicographical ordering allows one to compute the tilde 
vector polynomials along a strip of size $n+1$ in the $m$ direction. Equations 
\eqref{2.11} and \eqref{2.12} allow one to compute by changing $n$ and $m$ 
simultaneously. That is given the polynomials at the $(n,m-1)$ level and the 
$(n-1, m)$ level and the coefficients in \eqref{2.11}, \eqref{2.12} and 
\eqref{2.13} it is 
possible to compute $\PP_{n,m}$ and $\PPt_{n,m}$. It should be noted that 
in going from $(n-1,m)$, and $(n,m-1)$ to $(n,m)$ there are generically 4 new 
moments. Examination of the sizes of the coefficients in \eqref{2.11} and 
\eqref{2.12} shows that most of the coefficients in these matrices must be 
computable from the coefficients given at previous levels. Indeed in 
\cite{DGIM} it was shown that this is the case and that for $n\ge1$ and 
$m\ge1$ only $(\cK_{n,m})_{m,n}$, $(J^2_{n,m})_{m,n}$, $(J^1_{n,m})_{m,m}$ and 
$(J^1_{n,m})_{m,m+1}$ need be specified in order to compute $\PP_{n,m}$. For 
$n\ge1$ and $m\ge1$ we denote
\begin{align}\label{svalues}
(\cK_{n,m})_{m,n}=s_{2n-1,2m-1}\quad (J^2_{n,m})_{m,n}=s_{2n-1,2m}\\\
(J^1_{n,m})_{m,m}=s_{2n,2m-1}\quad(J^1_{n,m})_{m,m+1}=s_{2n,2m}.
\end{align}
For $n=0$ and $m\ge0$ or $n\ge1$ and $m=0$ we introduce $s_{0,m}$ and $s_{n,0}$
which are parameters in the one dimensional recurrence formulas associated 
with the line $(n,0)$ or $(0,m)$ respectively. Thus we find
\begin{Theorem}\label{th2.3}\cite{DGIM} Given parameters
$s_{0,0},\dots,s_{2n,2m}\in\R$, we construct
\begin{itemize}
\item  scalars $A_{i+1,0}$, $B_{i,0}$, $i=0,\dots,n-1$, and
$\At_{0,j+1},\Bt_{0,j}$, $j=0,\dots,m-1$;
\item  $j\times i$ matrices $\cK_{i,j}$ and $J^2_{i,j}$,
$i=1,\dots,n$, $j=1,\dots,m$;
\item $j\times (j+1)$ matrices $J^1_{i,j}$ for $i=1,\dots,n$,
$j=1,2\dots,m$.
\end{itemize}
If
\begin{equation}\label{2.17}
s_{2i,2j}>0 \text{ and }||\cK_{i,j}||<1,
\end{equation}
then there exists a positive linear functional $\cF$ such that
\begin{equation}\label{2.18}
\cF(\PP_{i,m},\PP_{j,m})=\delta_{i,j}I_{m+1}
\text{ and }
\cF(\PPt_{n,i},\PPt_{n,j})=\delta_{i,j}I_{n+1}.
\end{equation}
The conditions \eqref{2.17} are also necessary.
\end{Theorem}

With the above we are able to begin with the parameters $s_{i,j}$ and try to 
compute two variable orthogonality measures associated with these parameters. 
The connection with matrix orthogonal polynomials will be especially useful 
and we review some of the relevant relations. The theory of matrix orthogonal 
polynomials has an extensive literature \cite{Ber, Du, Ger, Mar, Ser, Si} and 
we begin with the equation~\eqref{2.10} and use \eqref{connectionform} to find,
\begin{equation}\label{matrixthreeterm}
  x P_{n,m} = A_{n+1,m} P_{n+1,m} + B_{n,m}P_{n,m} +
A_{n,m}^t P_{n-1,m}.
\end{equation}
If we set $P_{-1,m}=0$ and $P_{0,m}=I_{m+1}$ then the above equation uniquely 
gives $\{P_{n,m}\}$. We now turn the problem around and begin with 
$\{A_{n,m}\}_{n=1}^{\infty}$ and  $\{B_{n,m}\}_{n=0}^{\infty}$ where 
$A_{i,m}\ B_{i-1,m}$ are $(m+1)\times(m+1)$ real valued matrices with 
$B_{n,m}=B^t_{n,m}$ and $A_{n,m}$ lower triangular with strictly positive 
diagonal entries. We will suppose that
\begin{equation}\label{limitmatrix}
\lim_{n\to\infty} A_{n,m}=\frac{1}{2}I_{m+1}\ {\rm and} \lim_{n\to\infty} 
B_{n,m}=0,
\end{equation} 
and define $A_{0,m}=I_{m+1}$. In this case, there is a unique matrix measure 
$M$ such that 
\begin{equation}\label{orthonormal}
\int_{-1}^1 P_{n,m}(x)dM(x) P_{k,m}(x)^t =\delta_{n,k}I_{m+1}.
\end{equation}
Following \cite{Ger} we introduce the matrix valued function
\begin{equation}\label{psi}
\Psi_{n,m}(z)=P_{n,m}-2zA_{n,m}^tP_{n-1,m},
\end{equation}
where $z=x-\sqrt{x^2-1}$. The branch of the square root is chosen so that 
$z\rightarrow 0$ as $x\rightarrow +\infty$.
Using \eqref{matrixthreeterm} and \eqref{psi} one can deduce that 
\begin{equation}\label{2.23}
\begin{split}
&\Psi_{n,m}(z)\\&=\frac{1}{2z}A_{n,m}^{-1}\Psi_{n-1,m}+
\frac{1}{2}A_{n,m}^{-1}\left[(I_{m+1}-4A_{n,m}A_{n,m}^t)z-2B_{n-1,m}
\right]P_{n-1,m}.
\end{split}
\end{equation}
Note that if $A_{n,m}=\frac{1}{2} I_{m+1}$ and $B_{n-1,m}=0$ for all
$n\ge n_0$ then 
$$\Psi^*_{n,m}(z)=z^n\Psi_{n,m}(z)=\Psi^*_{n_0,m}(z) \text{ for all }
n\ge n_0.$$
We also introduce two scattering solutions of 
\eqref{matrixthreeterm} $P^{\pm}_{n,m}(z)$ that satisfy 
$$\lim_{n\to\infty} ||z^{\mp n}P^{\pm}_{n,m}(z)-I_{m+1}||=0$$
where the matrix norm is the Hilbert-Schmidt norm i.e 
$||B||=({\rm tr} (BB^{\dagger}))^{1/2}$. If we assume 
\begin{equation}\label{condone}
\sum_{n=1}^{\infty}n(||1-4A_{n,m}A_{n,m}^t||+||B_{n-1,m}||)<\infty,
\end{equation}
then a slight modification of the techniques in \cite{Ger} gives that for 
$n\ge0$, $P^+_{n,m}$ exists, is continuous for $|z|\le1$ and analytic for 
$|z|<1$ while $P^-_{n,m}$ exists, is continuous for $|z|\ge1$ and analytic for 
$|z|>1$. Since the coefficients in \eqref{matrixthreeterm} are real valued we 
see from the asymptotic conditions satisfied by $P^{\pm}_{n,m}$ 
that for $|z|=1$, $P^-_{n,m}(z)=P^+_{n,m}(1/z)$. For two solutions $Y$ and $X$ 
of \eqref{matrixthreeterm}
we define 
\begin{equation}\label{wronskian}
W_m[X,Y]=X_{n,m}^{\dagger}(\bar x)A_{n+1,m}Y_{n+1,m}(x)
-X_{n+1,m}^{\dagger}(\bar x)A^t_{n+1,m}Y_{n,m}(x),
\end{equation}
which by standard computations is independent of $n$. With the use of $W_m$ it 
follows that if \eqref{condone} holds then for all $|z|=1, z\ne\pm1$ 
\begin{equation}\label{expanp}
P_{n,m}(x)=\frac{2}{z-1/z}(P^+_{n,m}(z)f_-(z)- P^-_{n,m}(z)f_+(z))
\end{equation}
where 
\begin{equation}\label{jost}
f_+(z)=P^+_{-1,m}(z)^t \quad  |z|\le 1,
\end{equation}
and
\begin{equation}\label{jostpsi}
f_+(z)=\frac{1}{2z}\lim_{n\to\infty}\Psi^*_{n,m}(z),\ |z|\le 1.
\end{equation}
From the relation between $P^{\pm}_{n,m}$ we find,
\begin{equation}\label{ff}
f_+(z)=f_-(1/z)\quad |z|=1
\end{equation}
Manipulations similar to those that lead to \eqref{wronskian} yield
\begin{equation}
\begin{split}
&(P^+_{n,m}(z))^{\dagger}A_{n+1,m}P^+_{n+1,m}(z)
-(P^+_{n+1,m})^{\dagger}(z)A^t_{n+1,m}P^+_{n,m}(z)\\
&=\frac{1}{2}(z-1/z) I_{m+1}, 
\quad |z|=1, \quad z\ne \pm1,
\end{split}
\end{equation}
which implies that $P^+$ is nonsingular for $|z|=1, z\ne\pm1$. This follows 
since if there is a vector ${\bf a}$ and a $z_0, |z_0|=1, z_0\ne\pm1$ such 
that $P^+_{n,m}(z_0){\bf a}=0$ then $(P^+_{n,m}(z_0){\bf a})^{\dagger}=0$ 
which cannot happen by the equation above. This implies that $f_+$ is 
nonsingular for $|z|=1, z\ne\pm1$. With $n=-1$ in \eqref{expanp} we find for 
$|z|=1,\ z\ne\pm1$ the useful equation 
\begin{equation}\label{fpfm}
f^t_+ f_-=f^t_-f_+.
\end{equation} 
Another useful relation between $f_+$ and $P^+_{n,m}$ is 
\cite[formula 7.3]{Ger},
\begin{equation}\label{fp}
P^+_{n,m}(z)=\int_{-\infty}^{\infty}
   \frac{P_{n,m}(y)}{x-y}d M(y) f^t_+(z)\quad |z|<1
\end{equation}
where the relation between $x$ and $z$ is given above.

We will now make the following assumptions
\begin{equation}\label{assumone}
A_{n+1,m}=\frac{1}{2}I_{m+1},\quad B_{n,m}=0, \forall n\ge n_0,
\end{equation}
and
\begin{equation}\label{assumtwo}
\det(zf_+(z))\ne 0,\quad |z|\le 1,
\end{equation}
which will simplify the presentation below and are sufficiently general for 
our examples. The more general case will be taken up later (see however 
\cite{Ger} and \cite{Ser}). With assumption \eqref{assumone}
we see from \eqref{jostpsi} and \eqref{2.23} that 
$2zf_+(z)=\Psi^*_{n_0,m}(z)$ which will be a matrix polynomial of 
degree $2n_0$ if $A_{n_0,m}\ne\frac{1}{2}I_{m+1}$ and of degree $2n_0-1$ if 
$A_{n_0,m}=\frac{1}{2}I_{m+1}$ and $B_{n_0-1,m}\ne0$ and for $n\ge0$, 
$P^{\pm}_{n,m}$ are matrix polynomials in $z^{\pm 1}$ respectively. 
Furthermore \eqref{expanp} and contour integration show that the matrix 
orthogonal polynomials $\{P_{n,m}\}$ satisfy equation~\eqref{orthonormal} with
\begin{equation}\label{matrixweight}
dM(x)=\sigma_m (x) 
dx=\frac{1}{2\pi}\sqrt{1-x^2}(f_+(z)^{\dagger}f_+(z))^{-1} dx,
     \quad z=e^{i\theta}.
\end{equation}

What we would like to calculate next is what happens when we add a mass point 
to the system which can be accomplished by multiplying  $\Psi_{n_0,m}$ by a 
zero. To this end write,
\begin{equation}\label{psihat}
\hat\Psi_{n_0+1,m}(z)=(z_0/z-1)\Psi_{n_0,m}(z),
\end{equation} 
where $z_0$ is real and $|z_0|<1$.
With the choice of $\hat\Psi_{n_0+1,m}$ the computation of the orthogonality 
measure as well as the matrix orthogonal polynomials $\{\hat P_{n,m}(x)\}$ 
with $\hat P_{0,m}=I_{m+1}$ can be carried out in the following manner,
\begin{Theorem}\label{hat}
Given $\{A_{n,m}\}$ and $\{B_{n,m}\}$ satisfying \eqref{assumone} and 
$\Psi_{n_0,m}$ satisfying \eqref{assumtwo}. Let $\hat\Psi_{n_0+1,m}$ be given 
by \eqref{psihat}. Then $\{\hat P_{n,m}(x)\}$ are a set of matrix orthogonal 
polynomials satisfying
\begin{equation}\label{hatortho}
\int_{\R} \hat P_{n,m}(x)d\hat M(x) \hat P_{k,m}(x)^t dx
             =I_{m+1}\delta_{n,k} \quad n,k\ge0
\end{equation}
where
\begin{equation}\label{hatweight}
d\hat M(x)=\begin{cases}\frac{\hat\sigma_m(x)}{x_0-x},\quad -1<x<1\\
\hat r_m\delta(x-x_0),\quad x\notin[-1,1]\end{cases}
\end{equation}
with
\begin{equation}\label{hatsigma}
\hat\sigma_m(x)=\frac{\sigma_m(x)}{2 z_0}
\end{equation}
and 
\begin{equation}\label{hatr}
\hat r_m
=-\frac{(z_0-1/z_0)}{4z_0}[f_+(1/\bar z_0)^{\dagger}f_+(z_0)]^{-1}.
\end{equation}
The polynomials $\{\hat P_{n,m}\}$ have the Uvarov-Christoffel representation,
\begin{equation}\label{2.41}
\hat P_{n,m}(x)
=d_{n,m}(Q_{n-1,m}(x_0)^{\dagger}A_{n,m}P_{n,m}(x) 
          -Q_{n,m}(x_0)^{\dagger}A^t_{n,m}P_{n-1,m}(x)),
\end{equation}
where
\begin{equation}\label{q}
Q_{n,m}(x_0)=\int_{-\infty}^{\infty}P_{n,m}(x) d\hat M(x).
\end{equation}
The constants $d_{n,m}$ are chosen so that $\hat P_{n,m}$ is orthonormal for 
$n\ge0$, which is equivalent to
\begin{equation}\label{dnm}
(d_{n,m}^td_{n,m})^{-1}=\frac{1}{(2z_0)^3}
f_-(z_0)^{-1}P^-_{n-1,m}(z_0)^t A_{n,m}P^-_{n,m}(z_0)f^t_-(z_0)^{-1}.
\end{equation}

\end{Theorem}

\begin{proof} From the definition of $Q_{n,m}(x)$, $d\hat M(x)$ and 
$\hat\sigma_m$ we find 
for $n>0$ and $0<i<n$ that
\begin{equation}\label{ot}
\int_{-\infty}^{\infty}\hat P_{n,m}(x)(x-x_0)^i d\hat M(x)=0.
\end{equation}
Equations \eqref{fp}, \eqref{fpfm} and \eqref{expanp} show that
\begin{equation}\label{qform}
Q_n(x_0)=P^+_{n,m}(z_0)f_+^t(z_0)^{-1}/2z_0+P_{n,m}(x_0)\hat r_m
=P^-_{n,m}(z_0)f_-^t(z_0)^{-1}/2z_0,
\end{equation} 
where we have used the definition of $\hat r_m$ and the analytic properties 
of $P^{\pm}_{n,m}$. We use the above equation to define $Q_{-1,m}$. The 
substitution of this formula into \eqref{ot} with $i=0$ yields
\begin{equation}
\begin{split}
&d_{n,m}^{-1}\int_{-\infty}^{\infty}\hat P_{n,m}(x) d\hat M(x)\\
& =\frac{1}{(2 z_0)^2}(f_-(z_0)^{-1}(P^-_{n-1,m}(z_0)^t A_{n,m}P^+_{n,m}(z_0)
-P^-_{n,m}(z_0)^t A^t_{n,m}P^+_{n-1,m}(z_0))f_+^t(z_0)^{-1}\\
&+\frac{1}{2z_0}(f_-(z_0)^{-1}(P^-_{n-1,m}(z_0)^tA_{n,m}P_{n,m}(x_0)
-P^-_{n,m}(z_0)^tA^t_{n,m}P_{n-1,m}(x_0))\hat r_m.
\end{split}
\end{equation}
If in the first term on the left hand side we use \eqref{wronskian} and let 
$n$ tend to $\infty$ while in the second term we use \eqref{wronskian} and set 
$n=0$ we find from the definition of $\hat r_m$ and $f_+$ that the above 
integral is equal to 0 for $n>0$. To show that $d_{n,m}$ is given as in 
equation \eqref{dnm} we begin with the recurrence formula satisfied by 
$Q_{n,m}(x_0)$
\begin{equation}\label{qthreeterm}
 x_0 Q_{n,m} = A_{n+1,m} Q_{n+1,m} + B_{n,m}Q_{n,m} +
A_{n,m}^t Q_{n-1,m},
\end{equation}
which follows from \eqref{q}.
Thus routine manipulations give
$$
d^{-1}_{n,m}\hat P_{n,m}(x)=(x-x_0)
\sum_{i=0}^{n-1}Q^{\dagger}_{i,m}(x_0)P_{i,m}(x)+Q_{-1,m}(x_0),
$$
so that for $n>0$
\begin{equation}
\begin{split}
&\int_{\R}\hat P_{n,m}(x)d\hat M(x)\hat P^t_{n,m}(x)\\
&=-\frac{1}{2z_0}d_{n,m}\sum_{i=0}^{n-1}Q^{\dagger}_{i,m}(z_0)
\int_{-1}^{1}P_{i,m}(x)dM(x)\hat P^t_{n,m}(x)\\
&=\frac{1}{2z_0}d_{n,m}Q^t_{n-1,m}(z_0)A_{n,m}Q_{n,m}(z_0)d_{n,m}^t.
\end{split}
\end{equation}
The last equality follows by eliminating $\hat P_{n,m}$ using its definition.
\end{proof}

\begin{Remark}\label{re2.5} 
If we want to fix $d_{n,m}$ so that the highest coefficient of 
$\hat P_{n,m}(x)$ is a lower triangular matrix, we need to choose in 
\eqref{dnm} the unique solution for which $d_{n,m}Q_{n-1,m}(x_0)^{\dag}$ is 
a lower triangular matrix. Notice that for $n\geq n_0$ we have 
$P_{n,m}^{-}(z_0)=z_0^{-n}I_{m+1}$, $A_{n+1,m}=\frac{1}{2}I_{m+1}$ which leads 
to $d_{n+1,m}=4z_0^{n+2}f_-(z_0)$. Thus for $n>n_0$, formula \eqref{2.41} 
reduces to $\hat P_{n,m}=z_0P_{n,m}-P_{n-1,m}$. This implies 
$\hat A_{n+1,m}=\frac{1}{2}I_{m+1}$ and $\hat B_{n,m}=0$ for $n>n_0$.
\end{Remark}

\begin{Remark}\label{re2.6} Theorem \ref{hat} can be easily explained in 
terms of the Darboux transformation. Indeed, let $\cL_{m}$ be the second order 
difference operator corresponding to the second order difference equation 
\eqref{matrixthreeterm} for the polynomials $P_{n,m}(x)$, i.e. 
$(\cL_{m}(f))_n=A_{n+1,m} f_{n+1} + B_{n,m} f_{n} + A_{n,m}^t f_{n-1}$. 
We can think of $\cL_m$ as a second order difference operator acting on 
matrix valued functions (the size of the matrices is $(m+1)\times (m+1)$) 
of a discrete variable $n$. Similarly, we denote by $\hat{\cL}_m$ the 
operator corresponding to the polynomials $\hat P_{n,m}(x)$. 
Then using \eqref{qthreeterm} one can check that operator $\cL_m-x_0\Id$ can 
be factored as
\begin{equation}\label{2.49}
\cL_m-x_0\Id =\cP \cQ,
\end{equation}
where $\cQ$ is the backward difference operator in \eqref{2.41}, i.e.
\begin{equation}\label{2.50}
\cQ(f)_n=d_{n,m}\left(Q_{n-1,m}(x_0)^{\dagger}A_{n,m}f_n 
          -Q_{n,m}(x_0)^{\dagger}A^t_{n,m}f_{n-1}\right)
\end{equation}
and $\cP$ is the forward difference operator
\begin{equation*}
\cP(f)_n=(Q_{n,m}^{\dag})^{-1}
\left(f_{n+1}d_{n+1,m}^{-1}-f_{n}d_{n,m}^{-1}\right).
\end{equation*}
The operator $\hat{\cL}_m$ is obtained from $\cL_m$ by exchanging the 
factors in \eqref{2.49}, i.e. 
\begin{equation}\label{2.51}
\hat \cL_m-x_0\Id =\cQ \cP.
\end{equation}
\end{Remark}
The above discussion can be summarized as,
\begin{Theorem} Suppose that the coefficients in $\cL_m$ satisfy 
\eqref{assumone} and $zf_+(z)$ satisfies \eqref{assumtwo}. If $\hat\cL_m$ is 
related to $\cL_m$ by the Darboux transformation \eqref{2.49} and 
\eqref{2.51} and $\cQ$ is a backward difference operator of the form 
\eqref{2.50} where $Q_{n,m}(x_0)$ are given by \eqref{qform}
then the polynomials 
$\hat P_{n,m}$ are orthogonal with respect to $\hat M$
which is given by equations \eqref{hatortho}-\eqref{hatr}.
\end{Theorem} 

We can now use the above results to help solve the bivariate problem via the 
parametric moment problem. Equation \eqref{matrixmeasure} and 
Theorem~\ref{matrix} show that $dM^{m+1}$ is a Hankel matrix with entries 
$d\mu_j(x)=\int y^j d\mu(x,y)$ with the integration being over $y$. The above 
construction allows us to compute $\mu_j(x)$ for every $j$. 
We can think of $\mu_j(x)$ as one dimensional moments (in the 
variable $y$), depending on a parameter $x$. Let us consider the 
corresponding polynomials $q_m^x(y)$. The three term recurrence formula 
takes the form 
\begin{equation}\label{2.52}
yq^x_m(y)=a^x_{m+1}q^x_{m+1}(y)+b^x_{m}q^x_{m}(y)+a^x_{m}q^x_{m-1}(y),
\end{equation}
with coefficients depending also on the parameter $x$. 
Next we use one dimensional theory and compute $d\mu(x,y)$ by introducing 
the function 
\begin{equation}\label{2.53}
\psi^x_m(w)=q^x_m-2wa^x_mq^x_{m-1}.
\end{equation}
We will return to this general strategy in a later paper and will be content 
to illustrate it with the examples in the next section.

\section{One parameter deformation of the Chebyshev polynomials}

The two dimensional Chebyshev polynomials corresponding to the measure 
$$\frac{4}{\pi^2}\sqrt{1-x^2}\sqrt{1-y^2}dxdy$$ 
are parametrized in terms of the 
$s_{i,j}$ from Theorem~\ref{th2.3} as follows 
$s_{0,0}=1$, $s_{2n-1,0}=0$, $s_{2n,0}=\frac{1}{2}$, $s_{0,2m-1}=0$,  
$s_{0,2m}=\frac{1}{2}$, 
$s_{2n-1,2m-1}=0$, $s_{2n-1,2m}=0$, $s_{2n,2m-1}=0$, $s_{2n,2m}=\frac{1}{2}$ 
for $n,m\in\N$.

In order to have a nontrivial two dimensional polynomials (i.e. the measure is 
not just a product of two one dimensional measures) it was shown in 
\cite{DGIM} that  
$\cK_{i,j}\neq 0$ for at least one pair of indexes $(i,j)\in\N^2$. One way to 
construct such polynomials with relatively simple recurrence coefficients is 
as follows. We take $\cK_{1,1}=[s_{1,1}]$ to be a nonzero 
matrix, i.e. pick $s_{1,1}\neq 0$. Following the algorithm in Section 6 of 
\cite{DGIM}, we pick at each level $(n,m)$ with $n,m\in\N$ parameters 
$s_{2n-1,2m-1}=s_{1,1}$, $s_{2n-1,2m}=0$, $s_{2n,2m-1}=0$, 
$s_{2n,2m}=\frac{1}{2}$. 

This specific choice leads to a one parameter deformation of the Chebyshev 
polynomials with recurrence coefficients given by the following formulas
\begin{equation*}
\cK_{n,m}=\left[\begin{matrix}
0      & 0      &\cdots & 0 & 0\\
\vdots &\vdots  &    &\vdots& 0\\
0      & 0      & \cdots& 0& 0\\
0      & 0      & \cdots& 0 & s_{1,1}
\end{matrix}\right], 
\quad
J^{1}_{n,m}=\frac{1}{2}\left[\begin{matrix}
0& 1&  \\
1   & 0& 1 \\
        &         & \ddots & \ddots \\
        &         &       & 1 & 0 & 1 \\
        &          &      &        &  \sqrt{1-s_{1,1}^2} & 0 & 1
\end{matrix}\right],
\end{equation*}
and $J^2_{n,m}$ is a zero matrix. For $A_{n,m}$ and $B_{n,m}$ we obtain
\begin{equation*}
A_{1,m}=\frac{\sqrt{1-s_{1,1}^2}}{2}\left[\begin{matrix}
1 &   \\
        &         &\ddots & \\
        &         &       &  & 1 &  \\
        &          &      &        &   & \frac{1}{\sqrt{1-s_{1,1}^2}}
\end{matrix}\right],
\end{equation*}
\begin{equation*}
B_{0,m}= \frac{s_{1,1}}{2}\left[\begin{matrix}
0& 1&  \\
1   & 0& 1 \\
        &         & \ddots & \ddots\\
        &         &       & 1 & 0 & 1 \\
        &          &      &        &   1 & 0
\end{matrix}\right],
\end{equation*}
and for $n\geq 2$ we have 
\begin{equation}\label{3.1}
A_{n,m}=\frac{1}{2}I_{m+1}, \quad B_{n-1,m}=0,
\end{equation}
i.e. beyond $n=1$ the matrices $A_{n,m}$ and $B_{n-1,m}$ reach their 
asymptotic values. Similar formulas hold if we reverse the roles of 
$x$ and $y$, i.e. we have a measure symmetric in $x$ and $y$. 
 
Equation \eqref{3.1} shows that the right-hand side of formula \eqref{2.23} 
will vanish for $n\geq 2$, which gives that
\begin{equation*}
\Psi^*_{n,m}(z)=\Psi^*_{1,m}(z) \text{ for }n\geq 1,
\end{equation*}
and $f_+(z)=\frac{1}{2z}\Psi^*_{1,m}$ by \eqref{jostpsi}. This gives 
the matrix measure $dM(x)$ via formula \eqref{matrixweight}.

The entries of the  Hankel matrix $dM^{m+1}(x)$ allow us to compute 
the moments $\mu_j(x)dx=\int y^jd \mu_y(x,y)$, thus the coefficients 
in the three term recurrence relation \eqref{2.52}.

For the first few values one obtains that 
$$a^x_1=\frac{\sqrt{1-s_{1,1}^2}}{2}, \quad b^x_0=s_{1,1}x$$ 
and 
\begin{equation}\label{3.2}
a^x_{m}=\frac{1}{2}, \quad b^x_{m-1}=0, {\text{ for }}m\geq 2.
\end{equation}
This suggest one possible way to compute the measure
$d\mu(x,y)$, by first proving that \eqref{3.2} holds and then using 
one dimensional scattering theory in the variable $y$ with $x$ as a 
parameter. This would give that $w^m\psi^x_m(w)$ is independent of $m$ for 
$m\geq 1$, where $\psi^x_m(w)$ is the function defined by \eqref{2.53}, i.e.
$$w^m\psi^x_m(w)=w\psi^x_1(w)=
\frac{s_{1,1}^2w^2-2s_{1,1}xw+1}{\sqrt{1-s_{1,1}^2}}.$$
Applying the one dimensional scattering techniques, we can show that the 
measure is given by the following formula
\begin{equation}\label{3.3}
\begin{split}
d\mu(x,y)&=\frac{4}{\pi^2}\frac{\sqrt{1-x^2}\sqrt{1-y^2}}{\psi^x_m(w)\,
{\overline{{\psi^x_m(w)}}}}\,dxdy \\
&=\frac{4}{\pi^2}\frac{(1-s_{1,1}^2)\sqrt{1-x^2}\sqrt{1-y^2}}
{4s_{1,1}^2(x^2+y^2)-4s_{1,1}(1+s_{1,1}^2)xy+(1-s_{1,1}^2)^2}dx dy.
\end{split}
\end{equation}
In our particular case, we can give a simpler proof using the 
connection of the polynomials above with the Chebyshev polynomials. This is 
the content of the next theorem.

\begin{Theorem} \label{th3.1}
The polynomials $\PP_{n,m}(x,y)$ constructed with the parameters $s_{i,j}$ 
given by 
\begin{align*}
&s_{0,0}=1, \quad  s_{2n-1,0}=0, \quad s_{2n,0}=\frac{1}{2}, \quad 
s_{0,2m-1}=0, \quad  s_{0,2m}=\frac{1}{2}, \\
&s_{2n-1,2m-1}=s_{1,1}, \quad s_{2n-1,2m}=0, \quad s_{2n,2m-1}=0, 
\quad s_{2n,2m}=\frac{1}{2},
\end{align*}
for all $m,n\in\N$ are orthonormal with respect to the measure 
$d\mu(x,y)$ defined by \eqref{3.3}.
\end{Theorem}
\begin{proof}
We have
\begin{equation}\label{3.4}
\PP_{n,m}(x,y)=P_{n,m}(x)
\left[\begin{matrix} 1\\y\\ \vdots\\y^m
\end{matrix}\right].
\end{equation}
Thus the orthogonality relation 
\begin{equation}\label{3.5}
\int \PP_{n,m}(x,y)\PP_{k,m}^t(x,y)d\mu(x,y)= \delta_{n,k}I_{m+1}
\end{equation}
is equivalent to
\begin{equation}\label{3.6}
\int_{-1}^{1}\left( P_{n,m}(x)\int_{-1}^{1}
\left[\begin{matrix} 1\\y\\ \vdots\\y^m
\end{matrix}\right]\left[\begin{matrix} 1& y& \hdots &y^m
\end{matrix}\right] d_y\mu\right)P_{k,m}(x)^{t}\frac{2\sqrt{1-x^2}}{\pi}dx
=\delta_{n,k}I_{m+1},
\end{equation}
where 
$$d_y\mu=\frac{\pi}{2\sqrt{1-x^2}}d\mu$$ 
is a positive measure in the variable $y$. From \eqref{matrixweight}
we know that 
\begin{equation}\label{3.7}
\int_{-1}^{1}P_{n,m}(x) (\Psi_1^m(z)^{\dag}\Psi_1^m(z))^{-1} P_{k,m}(x)^{t} 
\frac{2\sqrt{1-x^2}}{\pi}dx =\delta_{n,k}I_{m+1}. 
\end{equation}
If we denote
\begin{equation}\label{3.8}
\phi_{n,m}(z,y)=\Psi_{n,m}(z)
\left[\begin{matrix} 1\\y\\ \vdots\\y^m
\end{matrix}\right],
\end{equation}
then to prove the orthogonality relation \eqref{3.5}, it is enough 
to show that 
\begin{equation}\label{3.9}
\int \phi_{1,m}(z,y)\phi_{1,m}(z,y)^{\dag} d_y\mu = I_{m+1}.
\end{equation}
The proof of the last equality can be easily established by using two facts. 
The first one is to notice that
$$\phi_{1,1}^0(z,y)=
\frac{s_{1,1}^2z^2-2s_{1,1}zy+1}{z\sqrt{1-s_{1,1}^2}},$$
i.e. for the measure given in \eqref{3.3} we have
\begin{equation}\label{3.10}
d_y\mu = \frac{2}{\pi}\frac{\sqrt{1-y^2}}{|\phi_{1,1}^0(z,y)|^2}dy.
\end{equation}
The second step is to connect $\phi_{1,m}(z,y)$ to the 
Chebyshev polynomials of the second kind. More precisely, we will show that 
\begin{equation}\label{3.11}
\phi_{1,m}^i(z,y)=
\left\{\begin{array}{ll} 
\phi_{1,1}^0(z,y)  U_i(y)         & \text{ if }i<m\\
\frac{1}{z}U_m(y)-s_{1,1}U_{m-1}(y)  &  \text{ if }i=m.
\end{array} \right.
\end{equation}
Notice that $1/\phi_{1,1}^0(z,y)$ is essentially a generating function for 
the polynomials $U_m(y)$. This combined with elementary properties of the 
Chebyshev polynomials leads to 
\begin{equation}\label{3.12}
\phi^m_{1,m}(z,y)=\sqrt{1-s_{1,1}^2}\phi_{1,1}^0(z,y)
\sum_{j=0}^{\infty}(zs_{1,1})^jU_{m+j}(y).
\end{equation}
Equation \eqref{3.9} follows immediately from \eqref{3.10}-\eqref{3.12}.

Thus it remains to establish \eqref{3.11}. From the recurrence coefficients, 
it is clear that 
\begin{equation}\label{3.13}
p_{0,m}^j(x,y)=\pt_{0,j}^0(x,y)=U_j(y) \text{ for every }j\in\N_0.
\end{equation}
From \eqref{2.11} with $n=1$ it follows that 
\begin{equation}\label{3.14}
p_{1,m}^j(x,y)=p_{1,m-1}^j(x,y) \text{ for } 0\leq j\leq m-2.
\end{equation}
This implies that 
\begin{equation}\label{3.15}
\phi_{1,m}^j(z,y)=\phi_{1,m-1}^j(z,y) \text{ for } 0\leq j\leq m-2.
\end{equation}
Since $J^3_{n,m}=-\cK_{n,m}\tilde A^t_{n-1,m}$ relation \eqref{2.12} for $n=1$ 
gives
\begin{subequations}
\begin{align}
&p_{1,m}^{j-1}+p_{1,m}^{j+1}=2yp_{1,m-1}^{j}=2yp_{1,m}^{j} \text{ for }
j\leq m-2 \label{3.16a}\\
&\sqrt{1-s_{1,1}^2}p_{1,m}^{m-2}+p_{1,m}^m=2yp_{1,m-1}^{m-1}
-s_{1,1}U_{m-1}(y), \label{3.16b}
\end{align}
\end{subequations}
where in the first equality we used \eqref{3.14} and in the second 
\eqref{3.13}.
Using the definition of $\phi_{n,m}$, \eqref{3.13} and \eqref{3.16a} 
it follows that
$$\phi_{1,m}^{j+1}+\phi_{1,m}^{j-1}=2y\phi_{1,m}^{j}\text{ for }j\leq m-2.$$
The last equation combined with the definition of $\phi_{1,m}$
leads by induction to $\phi_{1,m}^{j}=\phi_{1,1}^{0}U_j(y)$ for 
$j\leq m-1$ thus proving the first part in equation \eqref{3.11}. 

Finally, beginning with the definition of $\phi_{1,m}^m$ and then using 
equation \eqref{3.16b} to eliminate $p_{1,m}^m$ yields 
$$\phi_{1,m}^m=2yp_{1,m-1}^{m-1}-s_{1,1}U_{m-1}(y)
-\sqrt{1-s_{1,1}^2}p_{1,m}^{m-2}-zU_m(y).$$
Now substituting in the definition for $\phi^{m-1}_{1,m-1}$, 
using the recurrence formula for $2yU_{m-1}(y)$ and the definition for 
$\phi^{m-2}_{1,m}$ gives 
$$\phi_{1,m}^m=2y\phi_{1,m-1}^{m-1}+zs_{1,1}^2U_{m-2}(y)-s_{1,1}U_{m-1}(y)
-\sqrt{1-s_{1,1}^2}\phi_{1,m}^{m-2},$$
from which it follows
$$\phi_{1,m}^m-2y\phi_{1,m-1}^{m-1}=
\left(2ys_{1,1}-\frac{1}{z}\right)U_{m-2}(y)-s_{1,1}U_{m-1}(y),$$
which gives the second part in \eqref{3.11}.
\end{proof}

\begin{Remark}
We have computed also the first few coefficients for the recurrence relations 
\eqref{1.1}-\eqref{1.2} in the total degree ordering. The computations 
suggest the following formulas
\begin{equation}\label{3.17}
A_{x,n}=\frac{1}{2}\left[\begin{matrix}
s_{1,1} & \sqrt{1-s_{1,1}^2}&  \\
 &  & 1 \\
        &         &  & \ddots \\
        &         &       &  &   & 1 
\end{matrix}\right]
\end{equation}

\begin{equation}\label{3.18}
A_{y,n}=\frac{1}{2}\left[\begin{matrix}
1 &  &  \\
 & 1 &  \\
        &         &  &\ddots &\\
        &         &       &  &   1 & 0 
\end{matrix}\right],
\end{equation}
and $B_{x,n}$, $B_{y,n}$ are $(n+1)\times (n+1)$ zero matrices.
\end{Remark}

\section{Two parameter deformation of the Chebyshev polynomials}
Another deformation which leads to simple recurrence coefficients can 
be obtained by making the same choice of parameters as before except 
at levels $(1,0)$ and  $(0,1)$.
At level $(1,0)$ we leave $s_{1,0}$ free and we choose $s_{2,0}=\frac{1}{2}$.
At level $(0,1)$ we put $s_{0,1}=s_{1,0}s_{1,1}$ and $s_{0,2}=\frac{1}{2}$.
Using the same algorithm, we have this time two free 
parameters $s_{1,1}$ and $s_{1,0}$ and the recurrence coefficients take the 
form (we denote them by $A'$, $B'$, etc. in order to distinguish them from 
the one parameter deformation)
\begin{equation*}
A'_{1,m}=\frac{\sqrt{1-s_{1,1}^2}}{2}\left[\begin{matrix}
1 &   \\
-2s_{1,0}s_{1,1} &  1       & & \\
        &        &  \ddots    &   \ddots     &   & \\
    &&    &    -2s_{1,0}s_{1,1} &  1  \\
    &&&    &    -2s_{1,0}s_{1,1} &  \frac{1}{\sqrt{1-s_{1,1}^2}}  \\
\end{matrix}\right],
\end{equation*}
$A'_{n,m}=\frac{1}{2}I_{m+1}$ for $n\geq2$, 
\begin{equation*}
B'_{0,m}=\frac{1}{2}\left[\begin{matrix}
2s_{1,0} & s_{1,1}    \\
s_{1,1} &  2s_{1,0}(1-s_{1,1}^2)& s_{1,1}  \\
          &    &   \ddots     &    \\
& & s_{1,1} &  2s_{1,0}(1-s_{1,1}^2)& s_{1,1} \\
& & & s_{1,1} &  2s_{1,0}(1-s_{1,1}^2)     \\
\end{matrix}\right], 
\end{equation*}
 
\begin{equation*}
B'_{1,m}=\left[\begin{matrix}
s_{1,0}s_{1,1}^2 &   \\
 &  s_{1,0}s_{1,1}^2&  & \\
          &    &   \ddots     &    \\
& &  &  & s_{1,0}s_{1,1}^2   \\
& & &   &  & 0    \\
\end{matrix}\right], \quad B'_{n,m}=0 \text{ for }n\geq2,
\end{equation*} 
and the coefficients $\cK_{n,m}$ and $J^1_{n,m}$ are the same as in 
in the one parameter deformation. 
\begin{Theorem} \label{th4.1}
The polynomials $\PP'_{n,m}(x,y)$ corresponding to the parameters $s_{i,j}$ 
given by $s_{0,0}=1$, $|s_{1,0}|>\frac{1}{2}$, 
$s_{2,0}=\frac{1}{2}$,  $s_{2n-1,0}=0$, $s_{2n,0}=\frac{1}{2}$, for $n\geq 2$, 
$s_{0,1}=s_{1,0}s_{1,1}$, $s_{0,2}=\frac{1}{2}$, 
$s_{0,2m-1}=0$, $s_{0,2m}=\frac{1}{2}$ for $m\geq 2$, and 
\begin{equation*}
s_{2n-1,2m-1}=s_{1,1}, \quad s_{2n-1,2m}=0, \quad s_{2n,2m-1}=0, 
\quad s_{2n,2m}=\frac{1}{2}, \quad n,m\in\N
\end{equation*}
are orthonormal with respect to the measure 
\begin{equation}
\begin{split}
d\mu'(x,y)=&\frac{2z_0}{\pi^2(x_0-x)}\mu_0(x,y)\sqrt{1-x^2}\sqrt{1-y^2}dxdy\\
&\qquad +\frac{2(1-z_0^2)}{\pi}\mu_0(x,y)\delta(x-x_0)\sqrt{1-y^2}dx dy,
\end{split}
\end{equation}
where $z_0=1/(2s_{1,0})$ and 
\begin{equation}\label{4.2}
\mu_0(x,y)=
\frac{(1-s_{1,1}^2)}
{4s_{1,1}^2(x^2+y^2)-4s_{1,1}(1+s_{1,1}^2)xy+(1-s_{1,1}^2)^2}.
\end{equation}
\end{Theorem}

\begin{proof}
The shortest way to derive the orthogonality measure is to use 
Theorem~\ref{hat}. 
Indeed, let $\cL_{m}$ be the second order difference operator 
corresponding to the second order difference equation \eqref{2.10} for the 
polynomials $\PP_{n,m}(x,y)$ defined in Theorem~\ref{th3.1}, i.e. 
$(\cL_{m}(f))_n=A_{n+1,m} f_{n+1} + B_{n,m} f_{n} + A_{n,m}^t f_{n-1}$. 
Similarly, we denote by ${\cL'}_m$ the 
operator corresponding to the polynomials $\PP'_{n,m}(x,y)$ defined in 
Theorem \ref{th4.1}. Then 
one can check that $\cL_m$ and ${\cL'}_m$ are related by a Darboux 
transformation if we take 
$x_0=\frac{1}{2}(z_0+1/z_0)=(4s_{1,0}^2+1)/(4s_{1,0})$ and the 
operators $\cP$ and $\cQ$ are defined as follows
\begin{align*}
&\cP(f)_n=\frac{1}{2}f_{n+1}-\frac{1}{4s_{1,0}}f_n \text{ for }n\geq 1 
&&\text{ and }\cP(f)_0=A_{1,m}f_{1}-\frac{1}{4s_{1,0}}(\alpha^m)^t f_0 \\
&\cQ(f)_n=f_{n}-4s_{1,0}A^t_{n,m}f_{n-1} \text{ for }n\geq 1 
&&\text{ and } \cQ(f)_0=\alpha^m f_0,
\end{align*}
where $\alpha^m$ is the $(m+1)\times(m+1)$ matrix given by 
$$\alpha^m=\left[\begin{matrix}
1 &   \\
-2s_{1,0}s_{1,1} &  1       & & \\
 &-2s_{1,0}s_{1,1} &  1       & & \\
        &           &   \ddots     &   & \\
    &&&    &    -2s_{1,0}s_{1,1} &  1  \\
\end{matrix}\right].$$
From this relation it also follows that 
\begin{equation}\label{4.3}
\hat P_{n,m}(x)=z_0P'_{n,m}(x),
\end{equation} 
where $\hat P_{n,m}$ are the polynomials defined in Theorem \ref{hat}.
Writing the orthogonality relations for $\hat P_{n,m}(x)$, using \eqref{4.3}
and replacing $P'_{n,m}(x)[1\;y\dots y^m]^t$ by $\PP'_{n,m}(x,y)$
we obtained the desired orthogonality.
\end{proof}

\begin{Remark}
Finally, we list the recurrence coefficients in the total degree ordering. 
The matrices $A_{x,n}$ and 
$A_{y,n}$ are the same as the matrices for the 
polynomials in Section 3, i.e. they are given by \eqref{3.17} and 
\eqref{3.18}, respectively. For the $B_{x,n}$ we have 
$$B_{x.0}=[s_{1,0}], \quad 
B_{x,1}=
s_{1,0}\left[\begin{matrix}
1-s_{1,1}^2 & -s_{1,1}\sqrt{1-s_{1,1}^2}  \\
-s_{1,1}\sqrt{1-s_{1,1}^2} & s_{1,1}^2\\
\end{matrix}\right]$$
and for $n\geq 2$, $B_{x,n}$ is the block matrix
$$B_{x,n}=
\left[\begin{matrix}
B_{x,1} & 0 \\
0 & 0\\
\end{matrix}\right].$$
The matrices $B_{y,n}$ are identically equal to zero for $n\geq 1$ and 
$B_{y,0}=[s_{1,0}s_{1,1}]$.
\end{Remark}

\section*{Acknowledgments} 
We would like to thank A. Delgado, F. A. Gr\"unbaum and F. Marcell\'an for 
many useful discussions.

\end{document}